\documentclass[11pt]{article}
\usepackage{amssymb,amsmath,bm,graphicx,slashbox,array,multirow,rotating,subfigure}
\usepackage{enumerate}
\usepackage{amsthm}

\def\cN{{\mathcal N}}

\def\cF{{\mathcal F}}
\def\cL{{\mathcal L}}
\def\bR{{\mathbb R}}

\def\1{{\mathbf 1}}
\def\eps{{\mathbf \varepsilon}}
\def\E{{\mathbf E}}

\def\P{{\mathbf P}}

\def\arg{{\mathrm {arg}}}

\def\mle{{\mathrm{mle}}}
\def\b{{\mathrm {b}}}
\def\deq{{\stackrel{d}=}}
\def\wto{{\stackrel{w}\rightarrow}}
\def\Pto{{\stackrel{P}\rightarrow}}
\def\bL{{\mathbf L}}

\newtheorem{lemma}{Lemma}

\newtheorem{remark}{Remark}

\date{}

\title{On two estimates related to the change-point problem}
\author{\bfseries Farida Enikeeva}

\begin{document}
%\date{ }
\maketitle

\begin{center}
{\small\itshape $^1$The Kharkevich Institue for Information Transmission Problems of RAS,\\ 
Bol'shoi Karetnyi per. 19, GSP-4, 127994, Moscow, Russia}
\bigskip

\parbox{0.85\textwidth}{
\small
We consider the problem of estimating a smooth functional of an unknown signal with discontinuity from Gaussian observations. The signal is a known function that depends on an unknown parameter. This problem is closely related to the famous change-point problem. We obtain an asymptotic likelihood ratio process for the noise level tending to 0. Bayesian and maximum likelihood estimates are constructed and their relative efficiency is studied. Some simulation results and conclusions on non-asymptotic behavior of these estimates are presented.
\medskip

\medskip

}
\end{center}

\section{Introduction}

We consider the classical Gaussian white noise model 
\begin{equation}\label{wnoise}
dX_\eps(t) =f(t)dt +\eps dW(t),\quad t\in [0,1],
\end{equation}
where $f(t)$ is an unknown signal, $W(t)$ is a standard Brownian motion, the noise level $\eps>0$ is known. We assume that the function $f$ is continuous everywhere on $[0,1]$ except some unknown point $\tau$ and depends on some unknown parameter $\theta$, $f(t)\equiv f^\tau(\theta,t)$.

Let $\cL:\bL_2[0,1]\to \bR$ be a given smooth functional of $f$. The goal of this paper is to compare Bayesian and maximum likelihood estimates of $\cL[f]$ assuming that the function $f$ is known up to the parameters $\tau$ and $\theta$.
Let $\widehat\cL(X_\eps)$ be an estimate of $\cL[f]$. We will use the quadratic loss function and the mean squared risk for measuring the performance of the estimator:
$$
R_\eps(\widehat \cL,\cL)=\E_{\theta,\tau} (\widehat \cL(X_\eps)-\cL[f])^2.
$$

The model of observations~(\ref{wnoise}) of the Wiener process with a discontinuous drift was first considered by Ibragimov and Hasminski \cite{Ibragimov&Hasminski:1981}. Assuming that the function $f$ is known with an unknown discontinuity point $\tau$ the authors studied asymptotic efficiencies of Bayesian and maximum likelihood estimates of $\tau$ as $\eps\to 0$.  Asymptotic mean-square error of an MLE of the discontinuity point $\tau$ was calculated and an approximate value of the quadratic risk of a Bayes estimate was obtained. Later Rubin and Song \cite{Rubin&Song:1995} found an exact representation for the mean-square error of the Bayes estimate of $\tau$ in terms of Riemann's  zeta function. According to these results, the Bayes procedure of estimating a change point is asymptotically more efficient than the maximum likelihood procedure with the asymptotic relative efficiency~$\frac8{13}\zeta(3)\approx 0.7397$. 

This problem is closely related to the famous change-point problem considered by many authors. The literature on the change-point problem is vast, we refer the reader to the monographs of Cs\"org\H{o} and Horv\'ath \cite{Csorgo:1997} on asymptotic theory in the change-point problem, of Brodsky and Darhovsky \cite{Brodsky&Darhovsky:1993} on non-parametric methods, of Shiryaev \cite{Shiryaev:1978} on optimal detection of change in distribution, and many references therein. We also refer to an excellent review article of Bhattacharya \cite{Bhattacharya} that provides historical perspectives of the classical change-point problem. 

In spite of a long history of the change-point problem, the problem of estimating a smooth functional of a discontinuous signal was not considered.
We construct two estimates of $\cL$ for model~(\ref{wnoise}). We compare the asymptotic efficiencies of MLE and Bayesian estimate in the white noise model following the approach of \cite{Ibragimov&Hasminski:1981}.

The paper is organized as follows. In Section~2 we give a precise statement of the problem and obtain the asymptotic likelihood ratio process. In Section~3 the results on the relative efficiency of Bayesian and maximum likelihood estimates of the smooth functional are presented. Section~4 contains the results for a sequence version of~(\ref{wnoise}) with a simple signal representing the change in mean of a Gaussian sequence. In Section~5 we present simulation results for different signal-to-noise ratio and discuss both asymptotic and non-asymptotic aspects of the problem.

\section{Limiting likelihood ratio process}
It will be easier to work with a stochastic process $Y(t)$ satisfying the stochastic differential equation
\begin{equation}\label{wnoiseS}
dY(t)=\frac1\eps f^\tau(\theta,t)\,dt+\,dW(t),\quad t\in[0,1],
\end{equation}
where $W(t)$ is the standard Wiener process, $W(0)=0$, and $\eps>0$.

Assume that the functon $f^\tau(\theta,t)$ is defined as
\begin{equation}\label{signal}
f^\tau(\theta,t)=\begin{cases}
		f_1(\theta_1,t),& 0\le t<\tau\\
		f_2(\theta_2,t),& \tau< t\le 1,
	\end{cases}
\end{equation}
where $\theta=(\theta_1,\theta_2)$ and $\tau$ are unknown parameters that belong to some compact sets, $\theta\in\Theta=\Theta_1\times\Theta_2\subset \bR^2$, $\tau\in T=[a,b]$. We assume that $0<a\le\tau\le b<1$ so that that the change-point $\tau$ is separated from 0 and~1 and the change in the data happened within the interval $[a,b]$. 
 Denote by $\Delta=f^\tau(\theta,-\tau)-f^\tau(\theta,+\tau)\equiv f_1(\theta_1,\tau)-f_2(\theta_2,\tau)$ the jump size at the point $\tau$ assuming that $\Delta\neq0$.

In fact, $f^\tau(\theta,\cdot)$ depends on $\theta_1$ on $[0,\tau]$ and on $\theta_2$ on $[\tau,1]$. Thus, by abuse of notation we will write $\frac{\partial f_i}{\partial\theta}(\theta_i,t)$ meaning the value of the partial derivative of $f_i(x,t)$ with respect to $x$ at~$x=\theta_i$.

Let $\cF\subset\bL_2[0,1]$ be a linear space such that for any fixed parameters $\theta$ and $\tau$ the function $f^\tau(\theta,\cdot)\in\cF$ satisfies the following condition.

\noindent {\bf Condition F}. Assume that
\begin{enumerate}[(a)]
\item The functions $f_i(\theta_i,t)$, $i=1,2$, are continuous in $t$ on $[0,\tau]$ and $[\tau,1]$, and in $\theta_i$ on $\Theta_i$, respectively.
\item For any $x$ in a neighborhood of $\theta$, $f^\tau(x,t)$ has a bounded derivative $\frac{\partial f^\tau}{\partial t}(\theta,t)$ for all $t\in[0,1]$ except $t=\tau$.
%\item $f_i(\theta,t)$, $i=1,2$, is continuous and differentiable in $\theta$ at the point $\theta_i$ uniformly over $t$: $|\frac{\partial f_i}{\partial \theta}(\theta_i,t)|<\infty$ for all $t\in[0,1]$.
\item $f_i(\theta,t)$, $i=1,2$, are differentiable with respect to $\theta$ at $\theta_i$'s 
%such that $\sup\limits_{t\in[0,1]}|\frac{\partial f_i}{\partial \theta}(\theta_i,t)|<\infty$ and
such that
$$
\lim_{\delta\to 0} \frac1\delta\|f_i(\theta_i+\delta,\cdot)-f_i(\theta_i,\cdot)-\frac{\partial f_i}{\partial\theta}(\theta_i,\cdot)\delta\|_{\bL_2[0,1]}=0.
$$
\end{enumerate}

The problem is to estimate a smooth functional $\cL[f^\tau(\theta,\cdot)]$ of the signal $f^\tau(\theta,\cdot)$. Below the conditions on the functional $\cL$ are specified.

\noindent {\bf Condition L}. Let $\theta\in\Theta$ and $\tau\in T$ be fixed. The functional $\cL:\cF\subset\bL_2[0,1]\to \bR$ is Fr\'echet differentiable at $f^\tau(\theta,\cdot)\in\cF$.

\noindent {\bf Condition L$'$}. The value $\cL[f^\tau(\theta,\cdot)]$ of the functional $\cL:\cF\subset\bL_2[0,1]\to \bR$ at $f^\tau(\theta,\cdot)$ is differentiable with respect to $\tau$ in a neighborhood of $\theta$.

Let $\widehat\tau_\mle^\eps$ be an MLE of $\tau$ based on observations~(\ref{wnoiseS}). Let $\widehat\tau_\b^\eps$ be a Bayes estimate of $\tau$ based on observations~(\ref{wnoiseS}), where $\tau$ has some positive prior distribution on $[0,1]$. 
The analysis of quadratic errors of these two estimates is based on the properties of the stochastic process 
\begin{equation}\label{V}
V(t)=\exp(B(t)-|t|/2),
\end{equation}
where $B(t)$ is the two-sided Brownian motion defined by
\begin{equation}\label{two_sidedBM}
B(t)=\left\{\begin{array}{ll}
				W_1(t),& t\ge 0\\
				W_2(-t),& t<0.
				\end{array}
				\right.
\end{equation}
Here $W_i(t)$, $t\ge 0$, $i=1,2$ are independent standard Wiener processes  with $W_i(0)=0$. In fact, if $\Delta=f_1(\theta_1,\tau)-f_2(\theta_2,\tau)$, then the process $V(\Delta^2t)$ is a limiting process for the likelihood ratio of $\tau$ corresponding to the observations $Y(t)$ \cite{Ibragimov&Hasminski:1981}. 

Remind that $\theta$ and $\tau$ are defined on a compact  set $\Theta=\Theta_1\times\Theta_2\subset \bR^2$ and on the interval $T=(a,b)\subset [0,1]$, $0<a<b<1$, respectively. Following the approach of \cite{Ibragimov&Hasminski:1981} we will fix the unknown parameters $\theta=(\theta_1,\theta_2)$ and $\tau$ and  work with the local parameters  $h=(h_1,h_2)$ and $u$. Introduce the normalizing sets $\Theta_\eps=\Theta_\eps^1\times\Theta_\eps^2$, where $\Theta_\eps^i=\eps^{-1}(\Theta_i-\theta_i)$, $T_\eps=\eps^{-2} (T-\tau)$ such that $h\in \Theta_\eps$ and $u\in T_\eps$. 
Let $\P_{\theta,\tau}$ be the measure generated by the process~(\ref{wnoiseS}) and $Z^\eps_{\theta,\tau}(h,u)$ be the likelihood ratio of $\theta$ and $\tau$ based on this process,
$$
Z^\eps_{\theta,\tau}(h,u)= \frac{d\P_{\theta+\eps h,\tau+\eps^2u}}{d\P_{\theta,\tau}}(Y(t)).
$$

\begin{lemma}
Let $H=H_1\times H_2\subset \Theta_\eps$ and $U\subset T_\eps$ be compact sets and condition~F be satisfied. The distribution of the log-likelihood ratio process $\log \tilde Z^\eps_{\theta,\tau}(h,u)$ as $\eps\to 0$ converges uniformly over $(h,u)\in H\times U$ to the distribution of the process
\begin{equation}\label{logLL}
 \log  Z^0_{\theta,\tau} (h,u)= \frac12(Z_1^2+Z_2^2)-\frac12I_1^2\left(h_1-\frac{Z_1}{I_1}\right)^2\!\!\!-\frac12 I_2^2\left(h_2-\frac{Z_2}{I_2}\right)^2\!\!\!+\log V(\Delta^2 u),
\end{equation}
where $Z_1$ and $Z_2$ are independent $\cN(0,1)$, the process $V(u)$ defined in~(\ref{V}) is independent of $Z_1$ and $Z_2$, and
$$
I_1=\left(\int_0^\tau \left|\frac{\partial f_1}{\partial \theta}(\theta_1,t)\right|^2\,dt\right)^{1/2},\quad 
I_2=\left(\int_{\tau}^1 \left|\frac{\partial f_2}{\partial \theta}(\theta_2,t)\right|^2\,dt\right)^{1/2}.
$$
\end{lemma}
{\bf Proof}. 
From Girsanov's theorem (see~\cite{Ibragimov&Hasminski:1981}, Appendix~II, Theorem~1) it follows that the likelihood ratio for the measures generated by $Y(t)$ with the parameters $\theta$, $\tau$ and $\theta+\eps h$, $\tau+\eps^2 u$ satisfies
\begin{equation}\label{Girs}
\log Z_{\theta,\tau}^\eps(h,u)=\log\frac{d\P_{\theta+\eps h,\tau+\eps^2 u}}{d\P_{\theta,\tau}}(Y) =s_\eps(h,u)-\frac12r_\eps^2(h,u),
\end{equation}
where
\begin{gather*}
s_\eps(h,u)=\frac1\eps \int\limits_0^1 (f^{\tau+\eps^2 u}(\theta+\eps h,t)-f^\tau(\theta,t))\,d W(t),\\
r_\eps^2(h,u)=\frac1{\eps^2}\int_0^1(f^{\tau+\eps^2 u}(\theta+\eps h,t)-f^\tau(\theta,t))^2\,dt.
\end{gather*}
Using the same approach as in Lemma~7.2.1 of~\cite{Ibragimov&Hasminski:1981} and Condition~F[c] it is not difficult to show that uniformly over the compact  set $H\times U$ 
\begin{align*}
r_\eps^2(h,u)&=\Delta^2|u|+ h_1^2\int_0^\tau \left|\frac{\partial f_1}{\partial \theta}(\theta_1,t)\right|^2\,dt+h_2^2\int_{\tau}^1 \left|\frac{\partial f_2}{\partial \theta}(\theta_2,t)\right|^2\,dt+o(1),\\
&=\Delta^2|u|+ h_1^2 I_1^2+h_2^2I_2^2+o(1),
\quad\eps\to 0.
\end{align*}
Consider now the stochastic part $s_\eps(h,u)$ of the process. Let $u>0$. We have
\begin{align}
s_\eps(h,u))=
\frac1\eps \int\limits_0^\tau (f_1(\theta_1+\eps h_1,t)-f_1(\theta_1,t))\,dW(t)
&+\frac1\eps \int\limits_{\tau+\eps^2 u}^1 (f_2(\theta_2+\eps h_2,t)-f_2(\theta_2,t))\,dW(t)\nonumber\\
&+\frac1\eps \int\limits_\tau^{\tau+\eps^2 u} (f_1(\theta_1+\eps h_1,t)-f_2(\theta_2,t))\,dW(t)\label{seps}.
\end{align}

Following~\cite{Ibragimov&Hasminski:1981} from Condition~F[c] we can obtain the weak convergence of first two terms of~(\ref{seps})~to
$$
h_1Z_1 I_1\deq h_1\int\limits_0^\tau \frac{\partial f_1}{\partial\theta}(\theta_1,t)\,dW(t)\quad\mbox{and}\quad
h_2Z_2 I_2\deq h_2\int\limits_\tau^1 \frac{\partial f_2}{\partial\theta}(\theta_2,t)\,dW(t),
$$
respectively, where $I_1$ and $I_2$ are defined in the statement of the lemma and $Z_1$ and $Z_2$ are independent $\cN(0,1)$. 

Indeed, consider the first term of~(\ref{seps}). If
$$
V_\eps(h_1)=\frac1\eps \int\limits_0^\tau \left[(f_1(\theta_1+\eps h_1,t)-f_1(\theta_1,t))-h_1\eps \frac{\partial f_1}{\partial\theta}(\theta_1,t)\right]\,dW(t)
$$
then $\E V_\eps(h_1)=0$. From Condition~F[c] it follows that for some constant $C>0$ and any $h_1',h_1''\in H$
\begin{align*}
\E(V_\eps(h_1')-V_\eps(h_1''))^2\le C (h_1'-h_1'')^2 \int_0^\tau \left(\frac{\partial f_1}{\partial\theta}(\theta_1,t)\right)^2\,dt.
\end{align*}
%$
%\E(V_\eps(h_1')-V_\eps(h_1''))^2\le 2\sup\limits_{t\in[0,\tau]} \left(\frac{\partial f_1}{\partial\theta}(\theta_1,t)\right)^2 (h_1'-h_1'')^2.
%$
Consequently, from Theorem~1.A.19 of Prokhorov in~\cite{Ibragimov&Hasminski:1981}, p.~372 we obtain that for any $\delta>0$
$$
\lim_{\eps\to 0} \P_{\theta,\tau} \left\{\sup_{h_1\in H_1}|V_\eps(h_1)|>\delta\right\}=0.
$$
Next, from the properties of the stochastic integral we have
$$
\int\limits_0^\tau h_1\frac{\partial f_1}{\partial\theta}(\theta_1,t)\,dW(t)\deq h_1 Z_1\left(\int_0^\tau \left[\frac{\partial f_1}{\partial \theta}(\theta_1,t)\right]^2 dt\right)^{1/2}.
$$
Hence, uniformly over $h\in H$ the distribution of the first term in~(\ref{seps}) converges to the distribution of~$h_1Z_1 I_1$.
Similarly, we can show the weak convergence of the second term of~(\ref{seps}) to $h_2 Z_2 I_2$ uniformly over $H$.

Next, the last term in~(\ref{seps}) can be written as
\begin{align*}
\frac1\eps \int\limits_\tau^{\tau+\eps^2 u} (f_1(\theta_1+\eps h_1,t)-f_2(\theta_2,t))\,dW(t) &=
\frac1\eps \int\limits_\tau^{\tau+\eps^2 u} (f_1(\theta_1+\eps h_1,t)-f_1(\theta_1,\tau))\,dW(t)\\
&+\frac1\eps \int\limits_\tau^{\tau+\eps^2 u} (f_1(\theta_1,\tau)-f_2(\theta_2+\eps h_2,\tau))\,dW(t)\\
&+\frac1\eps \int\limits_\tau^{\tau+\eps^2 u} (f_2(\theta_2+\eps h_2,\tau)-f_2(\theta_2,t))\,dW(t).
\end{align*}
It can be shown, similarly to the proof in~\cite{Ibragimov&Hasminski:1981} that the first and the third terms converge to zero in probability uniformly over $H\times U$. For the second term we have
\begin{align*}
\frac1\eps \int\limits_\tau^{\tau+\eps^2 u} (f_1(\theta_1,\tau)-f_2(\theta_2+\eps h_2,\tau))\,dW(t)
&= \biggl(f_1(\theta_1,\tau)-f_2(\theta_2+\eps h_2,\tau)\biggr)\frac1\eps \bigl[W(\tau+\eps^2 u)-W(\tau)\bigr]\\
&\deq \biggl(f_1(\theta_1,\tau)-f_2(\theta_2+\eps h_2,\tau)\biggr) W_1(u) \Pto \Delta W_1(u).
\end{align*}
Note that the Wiener process $W_1(u)\deq \eps^{-1} \bigl[W(\tau+\eps^2 u)-W(\tau)\bigr]$ is independent of $Z_1$ and $Z_2$, since three summands in~(\ref{seps}) are independent.

Thus, combining the estimates for the stochastic and non-stochastic terms, we obtain the convergence of the distribution of $\log Z_{\theta,\tau}^\eps$ for $u>0$ to the distribution of 
\begin{align*}
\log Z^0_{\theta,\tau}(h,u)&= \Delta W_1(u) +Z_1 h_1 I_1+Z_2h_2 I_2-\frac12\Delta^2|u|-\frac12(h_1^2 I_1^2+h_2^2 I_2^2)\\
&=\frac12(Z_1^2+Z_2^2)-\frac12I_1^2\left(h_1-\frac{Z_1}{I_1}\right)^2\!\!\!-\frac12 I_2^2\left(h_2-\frac{Z_2}{I_2}\right)^2\!\!\!+\Delta\Bigl(B(u)-\frac12\Delta |u|\Bigr)
\end{align*}
uniformly over $H\times U$.
The similar analysis for $u\le 0$ yields the statement of the lemma.\hfill $\blacksquare$

\begin{remark}
This result can be generalized to the case of multiple change-points.
\end{remark}

\section{Relative efficiency of two estimates of $\cL[f]$}
First, recall the problem of estimating the point of discontinuity $\tau$ of the discontinuous signal $f$ from observations~(\ref{wnoise}) that was studied by Ibragimov and Hasminskii~\cite{Ibragimov&Hasminski:1981}. For quadratic loss function, they compared asymptotic efficiencies of a maximum likelihood and a generalized Bayesian estimators of $\tau$.  It turned out that  asymptotically the ratio of quadratic risks of Bayesian estimate and the MLE of $\tau$ does not depend on the function $f$ with the discontinuity point $\tau$ and that Bayesian estimator of $\tau$ is more efficient than the MLE of $\tau$. 

As it was mentioned above, if $\Delta$ is the jump size at the point $\tau$, then $V(\Delta^2 t)$ is a limiting likelihood ratio process for estimating $\tau$. 
Denote the MLE and Bayesian estimate of $\tau$ by $\widehat\tau_{\mle}^\eps$ and $\widehat\tau_\b^\eps$, respectively.
Let $u$ be a local parameter for $\tau$ with the normalization sequence $\eps^2$. Define
$$
\widehat\tau_{\mle}=\tau+\eps^2\widehat u_\mle,\quad \widehat\tau_\b=\tau+\eps^2\widehat u_\b,
$$
where  $\widehat u_\mle$ is the point at which the limiting likelihood attains its maximum and $\widehat u_\b$ is the generalized Bayesian estimate of $u$ that corresponds to the limiting likelihood. Namely, we have
\begin{equation}\label{u1u2}
\widehat u_\mle=\arg\max_{t\in\bR} V(\Delta^2 t)=\frac1{\Delta^2} \arg\max_{t\in\bR} V(t),\quad 
\widehat u_\b=\frac{\int_\bR t V(\Delta^2 t)\,dt}{\int_\bR V(\Delta^2 t)\,dt}=\frac1{\Delta^2}\frac{\int_\bR t V(t)\,dt}{\int_\bR V(t)\,dt}.
\end{equation}
Then the asymptotic relative efficiency of  $\widehat\tau_\mle^\eps$ and $\widehat\tau_\b^\eps$ coincides with the relative efficiency of the estimates $\widehat u_\mle$ and $\widehat u_\b$,
\begin{equation}\label{kappa0}
\lim_{\eps\to 0} \frac{\E_\tau(\widehat\tau_\b^\eps-\tau)^2}{\E_\tau (\widehat\tau_{\mle}^\eps-\tau)^2}=\lim_{\eps\to 0} \frac{\E_\tau(\widehat\tau_\b-\tau)^2}{\E_\tau (\widehat\tau_{\mle}-\tau)^2}=\frac{\E\,\widehat u_\b^2}{\E\,\widehat u_\mle^2}\equiv \kappa_0.
\end{equation}
Ibragimov and Hasminski~\cite{Ibragimov&Hasminski:1981} showed that $\E\,\widehat u_\mle^2=26/\Delta^4$, but they stated that $\E\,\widehat u_\b^2$ is hard to evaluate explicitly. Using computational methods, they obtained the following approximate value $\Delta^4\E\,\widehat  u_\b^2=19.5\pm 0.5$, and the efficiency
$
\kappa_0\approx 0.73\pm 0.03
$.
Later  Rubin and Song \cite{Rubin&Song:1995} obtained the exact value of $\E\widehat  u_\b^2$ and the asymptotic relative efficiency of two estimates,
$$
\kappa_0=\frac8{13}\zeta(3)\approx 0.7397,
$$
that appears to be very close to the approximate value found in~\cite{Ibragimov&Hasminski:1981}.  Here $\zeta$ is Riemann's zeta function defined as $\zeta(s)=\sum\limits_{n=1}^\infty n^{-s}$. 

We are interested in estimating the smooth functional $\cL[f^\tau(\theta,\cdot)]$. We will compare Bayesian and maximum likelihood estimates of $\cL[f^\tau(\theta,\cdot)]$. In fact, the problem is reduced to estimating the parameters $\theta=(\theta_1,\theta_2)$ and $\tau$.

We can write an MLE $\widehat \cL_\mle^\eps$ of $\cL[f^\tau(\theta,\cdot)]$ in terms of the local parameters $h$ and $u$ as 
$$
\widehat \cL_\mle^\eps=\cL[f^{\tau+\eps^2\widehat u^\eps}(\theta+\eps\widehat h^\eps,\cdot)]
$$
where
$$
(\widehat h^\eps,\widehat u^\eps)=\arg\max_{(h,u)\in\Theta_\eps\times T_\eps} Z^\eps_{\theta,\tau}(h,u).
$$
Define an MLE that corresponds to the limiting likelihood obtained in Lemma~1,
\begin{equation}\label{Lmle0}
\widehat \cL_\mle^0=\cL[f^{\tau+\eps^2\widehat u}(\theta+\eps\widehat h,\cdot)] 
\end{equation}
where $(\widehat h,\widehat u)\equiv(\widehat h_1,\widehat h_2,\widehat u)$ is the point at which the limiting log-likelihood process~(\ref{logLL}) attains its maximum. More precisely,
\begin{gather}
\widehat h_i=\arg\max_{h_i\in\bR} \left\{ \frac12{Z_i^2}-\frac12{I_i^2}\left(h_i-\frac{Z_i}{I_i}\right)^2\right\},\quad i=1,2,\label{hi}\\
\widehat u=\arg\max_{u\in\bR} V(\Delta^2 u)\equiv  \arg\max_{u\in\bR}\left\{ |\Delta| \left(B(u)-\frac{|\Delta u|}2\right)\right\}.\label{u}
\end{gather}

Obviously, $\E \widehat u=0$ and $\E \widehat h_i=0$,  $\E\widehat h_i^2=1/I_i^2$, $i=1,2$. We also know that $\E\widehat u^2=26/\Delta^4$.

Let $\widehat \cL_\b^\eps$  be a Bayesian estimate of $\cL[f^\tau(\theta,\cdot)]$ for quadratic loss function defined as
$$
\widehat \cL^\eps_{\b}=\arg\min_A \int\limits_{T_\eps}\int\limits_{\Theta_\eps}\bigl(A-\cL[f^{\tau+\eps^2u}(\theta+\eps h,\cdot)]\bigr)^2 Z^\eps_{\theta,\tau}(h,u)d h_1\,dh_2\,d u.
$$
Define the generalized Bayesian estimate  corresponding to the limiting likelihood process as
\begin{equation}\label{Lb0}
\widehat \cL^0_{\b}=\arg\min_{A} \int\limits_\bR \int\limits_{\bR^2}  \bigl(A- \cL[f^{\tau+\eps^2u}(\theta+\eps h,\cdot)] \bigr)^2 Z^0_{\theta,\tau} (h,u)\,dh_1\,dh_2\,du.
\end{equation}

Following the approach developed in \cite{Ibragimov&Hasminski:1981} we can prove the following result.

\begin{lemma}
Let $\theta\in\Theta$, $\tau\in T$, where $\Theta$ is a compact subset of $\bR^2$ and $T=[a,b]\subset[0,1]$. Let Conditions~F~and~L (L$'$) be satisfied for all $\theta\in\Theta$, $\tau\in T$. Then uniformly over $\Theta\times T$ 
\begin{equation}\label{rMLE}
\lim_{\eps\to 0}\frac{ \E_{\theta,\tau} \bigl(\widehat\cL_\mle^\eps-\cL[f^\tau(\theta,\cdot)]\bigr)^2}{\E_{\theta,\tau} \bigl(\widehat\cL_\mle^0-\cL[f^\tau(\theta,\cdot)]\bigr)^2}=1.
\end{equation}
and
\begin{equation}\label{rB}
\lim_{\eps\to 0} \frac{\E_{\theta,\tau} \bigl(\widehat\cL_\b^\eps-\cL[f^\tau(\theta,\cdot)]\bigr)^2}{\E_{\theta,\tau} \bigl(\widehat\cL_\b^0-\cL[f^\tau(\theta,\cdot)]\bigr)^2}=1.
\end{equation}
\end{lemma}
{\bf Proof}. The proof follows Theorems~1.10.1 and~1.10.2 of~\cite{Ibragimov&Hasminski:1981} and the continuity conditions~F and L (L$'$). We do not go into details here and just give a brief overview of the conditions of two theorems. Lemma~1 guarantees the convergence of the likelihood ratio on compact sets such that Condition~2 of Theorems~1.10.1 and~1.10.2 is satisfied. 

Next, by using of the technique developed in~\cite{Ibragimov&Hasminski:1981} it can be shown that for $\theta\in \Theta$, $\tau\in T=(a,b)\subset [0,1]$, $h',h''\in \Theta_\eps$, $u',u''\in T_\eps$, and $\eps<\eps_0$
\begin{align*}
\E_{\theta,\tau} \Bigl |[Z_{\theta,\tau}^\eps (h',u')]^{1/4}-[Z_{\theta,\tau}^\eps (h'',u'')]^{1/4}\Bigr|^4 
&\le \frac3{64} \frac1{\eps^2} \|f^{\tau+\eps^2 u'}(\theta+\eps h',\cdot)-f^{\tau+\eps^2 u''}(\theta+\eps h'',\cdot)\|^2_{\bL_2[0,1]}\\
&\le \frac3{64}\Bigl (C_1\|h'-h''\|^2+C_2|u'-u''|\Bigr)^2.
\end{align*}
It follows that for $\eps<\eps_0$ 
$$
\P_{\theta,\tau} \left\{Z_{\theta,\tau}^\eps (h,u)> \exp\biggl(-\frac{1}8 (C_1 \|h\|^2+C_2|\Delta||u|)\biggr)\right\}\le \exp\biggl(-\frac{1}{16}(C_1 \|h\|^2+C_2|\Delta||u|)\biggr).
$$
These two relations form Condition~1 of Theorems~1.10.1 and~1.10.2 that is necessary for consistency of the estimates. 

Finally, Condition~3 of the theorems is on the unique maximum of the limiting likelihood process~(\ref{logLL}) in case of an MLE of $\cL$ and on the unique minimum in~(\ref{Lb0}). From the results in~\cite{Ibragimov&Hasminski:1981} on the properites of the process $V(u)$ and concavity of  $\log \tilde Z^0_{\theta,\tau} (h,u)$ in $h$  it follows that these maximal (minimal) values are unique with probability~1.\hfill $\blacksquare$

\begin{lemma}
Let $\widehat \cL_{\mle}^\eps$ be an MLE and $\widehat \cL_{\b}^\eps$ be a generalized Bayesian estimate of $\cL[f^\tau(\theta,\cdot)]$, respectively. 
Assume that Conditions~L and~F are satisfied.
Then, asymptotically as $\eps\to 0$, the asymptotic quadratic risks of both estimators have the same first order term,
\begin{align*}
\lim_{\eps\to 0} \E_{\theta,\tau}\eps^{-2} \biggl(\widehat \cL_{\mle}^\eps-\cL[f^\tau(\theta,\cdot)]\biggr)^2
&=\lim_{\eps\to 0} \E_{\theta,\tau}\eps^{-2} \biggl(\widehat \cL_{\b}^\eps-\cL[f^\tau(\theta,\cdot)]\biggr)^2\\
&=\frac1{I_1^2} \biggl(\frac{\partial \cL}{\partial \theta_1}[f^\tau(\theta,\cdot)]\biggr)^2+\frac1{I_2^2} \biggl(\frac{\partial \cL}{\partial \theta_2}[f^\tau(\theta,\cdot)]\biggr)^2.
\end{align*}
\end{lemma}

{\it Proof}. 
From Lemma~2 it follows that  the risks~(\ref{rMLE}) and (\ref{rB}) of $\widehat\cL_{\mle}^\eps$ and  $\widehat\cL_{\b}^\eps$ have the same asymptotic behavior as the risks of $\widehat\cL_{\mle}^0$ and  $\widehat\cL_{\b}^0$. Thus we have to calculate the risks of the limiting estimates $\widehat \cL_{\mle}^0$  and $\widehat \cL_{\b}^0$ defined  in~(\ref{Lmle0}) and~(\ref{Lb0}), respectively. 

Remind that $\theta=(\theta_1,\theta_2)$ and $\tau$ are fixed and the increments $h=(h_1,h_2)$ and $u$ of $\theta$ and $\tau$ belong to some compact sets, $h\in H\subset \Theta_\eps$, $u\in U\subset T_\eps$. 

Let us first show that the following Taylor series expansion holds true
\begin{equation}\label{TaylorL}
\cL[f^{\tau+\eps^2 u}(\theta+\eps h,\cdot)] =\cL[f^\tau(\theta,\cdot)] + \frac{\partial \cL}{\partial\theta_1} [f^\tau(\theta,\cdot)]\eps h_1+ \frac{\partial \cL}{\partial\theta_2} [f^\tau(\theta,\cdot)]\eps h_2+o(\eps),\quad \eps\to 0.
\end{equation}

Since $\cL$ is Fr\'echet differentiable, for $f(t)\equiv f^\tau(\theta,t)\in\bL_2[0,1]$ there exists a linear mapping $\Lambda_{\theta,\tau}:\bL_2[0,1]\to \bR$ such that
\begin{equation}\label{1}
\cL[f^{\tau+\eps^2 u}(\theta+\eps h,\cdot)]-\cL[f^\tau(\theta,\cdot)]=\Lambda_{\theta,\tau}[\Delta f]+r[\Delta f]\|\Delta f\|_{\bL_2[0,1]},
\end{equation}
where 
$$
\Delta f(t)=f^{\tau+\eps^2 u}(\theta+\eps h,t)-f^\tau(\theta,t),\quad u\in U,\ h\in H
$$
and $\lim\limits_{\|\Delta f\|_{\bL_2[0,1]}\to 0} \|r(\Delta f)\|=0$ for the remainder term $r:\bL_2[0,1]\to \bR$.
Set
$$
\Delta f(t)=\bigl(f^{\tau+\eps^2 u}(\theta+\eps h,t)-f^{\tau}(\theta+\eps h,t)\bigr)+\bigl(f^{\tau}(\theta+\eps h,t)-f^\tau(\theta,t)\bigr)\equiv \Delta_\theta f(t)+\Delta_\tau f(t).
$$

Since $H$ and $U$ are compact sets and $f$ is continuous in $t$ everywhere except $t=\tau$ (Conditions~F[a-b]), we have 
$$
\Delta_\theta f(t)=f^{\tau+\eps^2 u}(\theta+\eps h,t)-f^{\tau}(\theta+\eps h,t)=r_0^\eps(\theta,\tau,h,u,t)
$$
where $\lim\limits_{\eps\to 0} \sup\limits_{h\in H,\ u\in U}\|r_0^\eps(\theta,\tau,h,u,\cdot)\|_{\bL_2[0,1]}=0$. Indeed, we have for $t\neq \tau$
$$
\Delta_\theta f(t)=(f_1(\theta_1+\eps h,t)-f_2(\theta_2+\eps h_2,t))\1\{|t-\tau|<\eps^2 u\}(\1\{u>0\}-\1\{u<0\}),
$$
where the functions $f_1$ and $f_2$ are continuous on $[0,1]$. 

From Condition~F[c] it follows that  for $t\neq \tau$
\begin{align*}
\Delta_\tau f(t)&=f^{\tau}(\theta_1+\eps h_1,\theta_2+\eps h_2,t)-f^\tau(\theta_1,\theta_2,t)\\
&=\frac{\partial f^{\tau}(\theta,t)}{\partial\theta_1}  \eps h_1+\frac{\partial  f^{\tau}(\theta,t)}{\partial\theta_2} \eps h_2+r_1(\eps h_1,t)\eps |h_1|\1\{t<\tau\}
+r_2(\eps h_2,t)\eps |h_2|\1\{t>\tau\}
\end{align*}
where  $\lim\limits_{\eps\to 0}\sup\limits_{h_i\in H_i} \|r_i(\eps h_i,\cdot)\|_{\bL_2[0,1]}=0$.
Combining two formulas for $\Delta_\tau f(t)$ and $\Delta_\theta f(t)$ we obtain
$$
\Delta f(t)=\frac{\partial f^{\tau}(\theta_1,t)}{\partial\theta_1}  \eps h_1+\frac{\partial  f^{\tau}(\theta_2,t)}{\partial\theta_2} \eps h_2+r^\eps(\theta,\tau,h,u,t)
$$
where $\sup\limits_{h\in H,\ u\in U}\|r^\eps(\theta,\tau,h,u,\cdot)\|_{\bL_2[0,1]}\to 0$, $\eps\to 0$. Thus, $\|\Delta f\|_{\bL_2[0,1]}\le \|\Delta_\theta f\|_{\bL_2[0,1]}+\|\Delta_\tau f\|_{\bL_2[0,1]}\to 0$ as $\eps \to 0$ uniformly over $h$ and $u$ and, consequently, $\|r[\Delta f]\|=o(\eps)$ in~(\ref{1}).

Substituting the obtained expansion in~(\ref{1}) gives the desired formula~(\ref{TaylorL}), where  the linear mappings are defined as $\frac{\partial \cL}{\partial \theta_i}[f^\tau(\theta,\cdot)]=\frac{\partial f^{\tau}(\theta,\cdot)}{\partial\theta_i}\circ \Lambda_{\theta,\tau}: \Theta_i\to \bR$.

First, we will calculate the  asymptotic mean-square error of $\widehat \cL_{\mle}^0$ defined in~(\ref{Lmle0}). From formula~(\ref{TaylorL}) and independence of $\widehat h_1$ and $\widehat h_2$ we obtain
\begin{align*}
\E_{\theta,\tau} (\widehat \cL_{\mle}^0-\cL[f^\tau(\theta,\cdot)])^2&=\eps^2 \biggl(\frac{\partial \cL}{\partial\theta_1} [f^\tau(\theta,\cdot)]\biggr)^2 \E_{\theta,\tau} \widehat h_1^2 +\eps^2 \biggl(\frac{\partial \cL}{\partial\theta_2} [f^\tau(\theta,\cdot)] \biggr)^2\E_{\theta,\tau} \widehat h_1^2 +o(\eps^2)\\
&= \frac{\eps^2}{I_1^2} \biggl(\frac{\partial \cL}{\partial \theta_1}[f^\tau(\theta,\cdot)]\biggr)^2+\frac{\eps^2}{I_2^2} \biggl(\frac{\partial \cL}{\partial \theta_2}[f^\tau(\theta,\cdot)]\biggr)^2+o(\eps^2),
\end{align*}
since $\E_{\theta,\tau} \widehat h_i^2=1/I_i^2$ and $\E_{\theta,\tau} \widehat h_i=0$.

Now calculate the risk of Bayesian estimate $\widehat \cL_{\b}^0$ defined in~(\ref{Lb0}). Using the Taylor series expansion~(\ref{TaylorL}) and formula~(\ref{logLL}) we have
\begin{align*}
\widehat \cL_{\b}^0&=\frac{\displaystyle{\int\limits_\bR\int\limits_{\bR^2} \cL[f^{\tau+\eps^2u}(\theta+\eps h,\cdot)] \tilde Z^0_{\theta,\tau} (h,u)\,dh_1\,dh_2\,du}}
{\displaystyle{\int\limits_\bR\int\limits_{\bR^2} \tilde Z^0_{\theta,\tau} (h,u)\,dh_1\,dh_2\,du}}\\
&\qquad\qquad\qquad=\cL[f^\tau(\theta,\cdot)]+\eps\frac{\partial \cL}{\partial\theta_1}[f^\tau(\theta,\cdot)] \frac{Z_1}{I_1} +\eps\frac{\partial \cL}{\partial\theta_2}[f^\tau(\theta,\cdot)] \frac{Z_2}{I_2}+o(\eps),
\end{align*}
where $Z_i$'s are independent $\cN(0,1)$. Calculating the mean square risk gives exactly the same asymptotic behavior as the one for the risk of $\widehat \cL_\mle^0$. \hfill $\blacksquare$

\begin{remark}
If Condition~L$'$ is satisfied and $\cL$ is twice differentiable w.r.t.\ $\theta_i$, we can calculate the second order terms of the asymptotic risks applying the facts that $\E \widehat u_\b=\E\widehat u_\mle=0$, $\E\widehat u_\mle^2=26/\Delta^4$, $\E\widehat u_\b^2=16\zeta(3)/\Delta^4$. We have
\begin{align*}
\E_{\theta,\tau} (\widehat \cL_{\mle}^\eps-\cL[f^\tau(\theta,\cdot)])^2=\eps^2\sum_{i=1,2}\frac1{I_i^2} \biggl(\frac{\partial \cL}{\partial \theta_i}[f^\tau(\theta,\cdot)]\biggr)^2
&+\eps^4 \frac{26}{\Delta^4} \biggl(\frac{\partial \cL}{\partial \tau}[f^\tau(\theta,\cdot)]\biggr)^2\\ &+\eps^4 \sum_{i=1,2}\frac3{I_i^4} \biggl(\frac{\partial^2 \cL}{\partial\theta_i^2}[f^\tau(\theta,\cdot)]\biggr)^2+o(\eps^4),\\
\E_{\theta,\tau} (\widehat \cL_{\b}^\eps-\cL[f^\tau(\theta,\cdot)])^2=\eps^2\sum_{i=1,2}\frac1{I_i^2} \biggl(\frac{\partial \cL}{\partial \theta_i}[f^\tau(\theta,\cdot)]\biggr)^2
&+\eps^4 \frac{16\zeta(3)}{\Delta^4} \biggl(\frac{\partial \cL}{\partial \tau}[f^\tau(\theta,\cdot)]\biggr)^2\\ &+\eps^4 \sum_{i=1,2}\frac{I_i^4+2I_i^2+3}{I_i^4} \biggl(\frac{\partial^2 \cL}{\partial\theta_i^2}[f^\tau(\theta,\cdot)]\biggr)^2+o(\eps^4).
\end{align*}
If Condition~L$'$ is satisfied and $\cL[f^\tau(\theta,\cdot)]$ is a function of $\tau$ only, $\cL[f^\tau(\theta,\cdot)]\equiv g(\tau)$,  we obtain the result~(\ref{kappa0})  of Ibragimov and Hasminski~as a corollary:
$$
\lim_{\eps\to 0}\frac{\E_{\theta,\tau} (\widehat \cL_{\b}^\eps-\cL[f^\tau(\theta,\cdot)])^2}{\E_{\theta,\tau} (\widehat \cL_{\mle}^\eps-\cL[f^\tau(\theta,\cdot)])^2}=
\lim_{\eps\to 0}\frac{\E_{\theta,\tau} (\widehat \cL_{\b}^\eps-g(\tau))^2}{\E_{\theta,\tau} (\widehat \cL_{\mle}^\eps-g(\tau))^2}=\frac8{13}\zeta(3).
$$
\end{remark}

\section{Estimation in the sequence model}
In this section we give explicit estimates for the problem of estimating a smooth functional in the equivalent sequence model. We assume that a very simple signal is observed, which is constant up to some moment of time $\tau$ and equals zero afterwards.

We observe the vector $X=(X_1,\dots,X_n)$, where $X_i$'s are Gaussian random variables with distribution $\P_{\theta_i,\tau}$ defined on the probability space $(\mathcal X,\mathcal B_i, \P_{\theta_i,\tau})$,
\begin{equation}\label{model_chp}
X_i=\theta_i+\eps\xi_i,\quad i=1,\dots, N.
\end{equation}
The signal is constant up to some moment of time $\tau$, $\theta_i=\theta\1(i-\tau\le 0)$. The parameters $\tau\in\{1,\dots,n\}$ and $\theta$ are unknown, $\xi_i$'s are i.i.d. $\cN(0,1)$, $\eps>0$ is known. The goal is to estimate a smooth function $L(\theta,\tau)$ of the signal. 

\subsection{Overview}

The change-point problem for the change in mean in Bayesian set-up was considered by Chernoff and Zacks \cite{Chernoff&Zacks:1964}. They studied the change in mean for a sequence of Gaussian r.v.'s and obtained a Bayes estimate of the difference in means before and after the change. This estimate was obtained under the assumption that the current mean and the jump size have normal prior distributions. The change-point was also assumed to be random with some arbitrary discrete prior.  In their paper, Chernoff and Zacks compare Bayesian estimates of the current mean and the minimum variance linear unbiased estimates (MVLUE) when the signal-to-noise ratio is greater than~2. In particular, they state that if there is exactly one change in the observations, the Bayes estimator is very efficient if the change takes place at the beginning of the sequence and it looses its efficiency when the change takes place very close to the last observation. In the case of at most one change  Bayesian procedures are always better than MVLUE. Rubin  \cite{Rubin:1961} considered the change-point problem in the context of estimating discontinuities in multivariate densities. He discussed Bayesian and maximum likelihood approach to the problem and mentioned that the Bayes procedure with the uniform prior distribution on the unknown parameter $\tau$ gives more efficient estimates than the maximum likelihood approach.

A maximum likelihood estimate (MLE) of a change-point was first obtained by Hinkley \cite{Hinkley:1970} in the problem of estimating a moment of the change in mean of Gaussian data under the assumption that the jump size  is small.  Later an asymptotic distribution of the MLE of a change-point for the case of close normal means was derived by Bhattacharya and Brockwell \cite{Bhattacharya&Brockwell:1976}. Their result was generalized by Bhattacharya \cite{Bhattacharya:1987} to the case of a small jump size in a multidimensional parameter.  Following \cite{Bhattacharya&Brockwell:1976} and \cite{Bhattacharya:1987} Ferger \cite{Ferger:1994} proposed a class of estimates  for a change-point based on $U$-statistics for the case of small disorders in the distribution. 
Later Brodskii and Darhovskii \cite{Brodskii&Darhovskii:1990} studied an asymptotic behavior of an estimate for the change-point in Gaussian sequence with unknown mean without the assumption that the difference between the means (the jump size) tends to zero. They proposed a family of estimates for the change-point based on the Kolmogorov--Smirnov statistics which includes an MLE.  An asymptotic distribution of these estimates was derived and the corresponding testing problem was considered.

Let $\nu^n$ be a $\sigma$-finite measure on $\sigma$-algebra $\mathcal B=\mathcal B_1\times\dots\times\mathcal B_n$ and $\P_{\theta,\tau}^n=\P_{\theta_1,\tau}\times\dots\times\P_{\theta_n,\tau}$. Then the joint density of $X$ (likelihood) is given by
\begin{equation}\label{density}
\frac{d\P_{\theta,\tau}^n}{d\nu^n}(X)\equiv p_n^\eps(X;\theta,\tau)=(2\pi\eps^2)^{-n/2} \exp\left\{-\frac1{2\eps^2} \left(\sum\limits_{i=1}^\tau (X_i-\theta)^2+\sum\limits_{i=\tau+1}^n X_i^2\right)\right\}.
\end{equation}

Assuming that $\eps=1$ in model~(\ref{model_chp}), Bhattacharya and Brockwell \cite{Bhattacharya&Brockwell:1976} derived a limiting process for the likelihood ratio under the conditions that the parameter $\theta$ is small, $\theta=\delta\nu_n^{-1}$, where $\nu_n\to\infty$ slower than $n^{1/2}$ and the length of the observed sequence $n\to\infty$. According to their result, as $n\to\infty$, the following weak convergence holds with respect to uniform convergence on compact sets,
$$
\log\frac{d\P_{\theta+n^{-1/2} h,\tau+\nu_n^2 u}}{d\P_{\theta,\tau}}
\wto \frac{Z^2}2 -\frac 12\lambda \left( h-\frac{Z}{\sqrt\lambda}\right)^2 +|\delta |\left(B(u)-\frac12 |\delta| |u|\right),\quad n\to\infty
$$
where $Z$ is $\cN(0,1)$ and $B(u)$ is a two-sided Wiener process~(\ref{two_sidedBM}) independent of $Z$ and $\lambda=\lim\limits_{n\to\infty}\tau/n$.

\begin{remark}
Some general results on the behavior of the log-likelihood ratio and of the maximum likelihood estimate of the change point $\tau$ can be found in~Section~1.6 of~\cite{Csorgo:1997}. In particular, Theorems~1.6.2 and~1.6.3 of~\cite{Csorgo:1997} state that the asymptotic distribution as $n\to\infty$ of the likelihood ratio in the situation of a decreasing size of the change in means obtained in \cite{Bhattacharya&Brockwell:1976} differs from the one in the situation of a fixed change in mean.
\end{remark}

In our case the number of observations $n$ is fixed and the change in mean (the jump size at $\tau$) $\theta$ is fixed. Let $\theta\in\Theta$, where $\Theta$ is an open subset of $\bR$ and $\tau=[n\alpha]$, where $\alpha\in A=(0,1)$.  Define the sets $\Theta_\eps=\eps^{-1}(\Theta-\theta)$ and $A_\eps=\eps^{-2}(A-\alpha)$. The following lemma that is given without proof describes the asymptotic behavior of the likelihood ratio as $\eps\to 0$.

\begin{lemma}
Let $X=(X_1,\dots,X_n)$ be given by~(\ref{model_chp}) and $H$, $V$ be some compact subsets of $\Theta_\eps$ and $A_\eps$, respectively. Then the following weak convergence holds uniformly over $(h,v)\in H\times V$ as $\eps\to 0$,
\begin{equation}\label{limitZ}
\log \frac{d\P_{\theta+\eps h,\tau+\eps^2 nv}}{d\P_{\theta,\tau}}(X)\wto \frac{Z^2}2-\frac\tau2 \left(h-\frac{Z}{\sqrt\tau}\right)^2+|\theta|\sqrt n\left( B(v)-\frac12 \sqrt n|\theta| |v|\right),\quad \eps\to 0.
\end{equation}
where $Z$ is $\cN(0,1)$ independent of the two-sided Brownian motion $B(v)$.
\end{lemma}

\subsection{MLE and Bayesian estimate of $L(\theta,\tau)$}
Let us find an MLE of $L(\theta,\tau)$. The log-likelihood $\log p_n^\eps(X;\theta,\tau)$ satisfies
\begin{eqnarray*}
\log p_n^\eps(X;\theta,\tau)+\frac n2 \log(2\pi\eps^2)&=&-\sum\limits_{i=1}^\tau\frac{(X_i-\theta)^2}{2\eps^2}-\sum\limits_{i=\tau+1}^n\frac{X_i^2}{2\eps^2}\\
&=&-\frac1{2\eps^2}\sum\limits_{i=1}^n X_i^2 +\frac\theta{\eps^2}\sum\limits_{i=1}^\tau X_i-\frac{\theta^2\tau}{2\eps^2}.
\end{eqnarray*}
First, we maximize the log-likelihood with respect to $\theta$ and replace $\theta$ by its conditional MLE $\bar X_\tau=\frac1\tau\sum_{i=1}^\tau X_i$. 
Next, maximizing the obtained log-likelihood with respect to $\tau$  we obtain an MLE of the change-point $\tau$,
\begin{equation}\label{MLEtau}
\widehat\tau_{\mathrm{mle}}^\eps=\arg\max_{1\le k\le n} \left\{ \frac1{2\eps^2 k} \left(\sum\limits_{i=1}^k X_i\right)^2\right\} =\arg\max_{1\le k\le n} U_k,
\end{equation}
where
\begin{equation}\label{Uk}
U_k=\frac1{2\eps^2 k} \left(\sum\limits_{i=1}^k X_i\right)^2.
\end{equation}
The estimate~(\ref{MLEtau}) of $\tau$ for unknown  change in mean $\theta$ of normal distribution  was first obtained by Hinkley in \cite{Hinkley:1970}. 
An asymptotic distribution as $\eps\to 0$ of this estimate was derived by Brodskii and Darkhovskii in~\cite{Brodskii&Darhovskii:1990}. 

Finally, an MLE of $L(\theta,\tau)$ is given by
$$
\widehat L_{\mle}^\eps=L(\widehat\theta_{\mle},\widehat\tau_{\mle})=L(\bar X_{\widehat\tau_{\mle}},\widehat\tau_{\mle}).
$$
For example, if $L(\theta,\tau)=\sum_{i=1}^\tau\theta_i=\theta\tau$, then
$
\widehat L_{\mle}^\eps=\sum\limits_{i=1}^{\widehat\tau_{\mathrm{mle}}} X_i
$.

\begin{remark}
Note that using~(\ref{model_chp}) we can write
$$
\widehat \tau_\mle=\arg\max_{1\le k\le n}\left\{ \Bigl| \frac{\eps}{\sqrt k}\sum_{i=1}^{k} \xi_i+\theta \sqrt k\1\{k\le\tau\}+\frac{\theta\tau}{\sqrt k}\1\{k>\tau\}\Bigr|\right\}.
$$
Thus, to find a non-asymptotic risk of $\widehat L_\mle$ we need to calculate the joint distribution of $\widehat\tau_\mle$ and $\sum_{i=1}^{\widehat\tau_\mle} \xi_i$. This problem is similar to calculation of the joint distribution of
$$
\widehat \tau=\arg\max\limits_{t>0} \Bigl|\theta\sqrt t\1\{t\le\tau\}+\frac{\theta \tau}{\sqrt t}\1\{t>\tau\}+\eps \frac{W(t)}{\sqrt t}\Bigr|\quad\mbox{and $W(\widehat\tau)$.}
$$ 
\end{remark}

Let us now find a Bayesian estimate of $L(\theta,\tau)$. Assume that $\theta$ has a non-informative prior distribution $\cN(0,\sigma^2)$, where $\sigma^2\to\infty$ and $\tau$ is uniformly distributed on the set $\{1,\dots,n\}$. Then the generalized posterior density of $(\theta,\tau)$ is given by
\begin{equation}\label{posterior}
\pi(\theta,\tau|X)=\frac{1}{\sqrt{2\pi\eps^2}} \frac{\displaystyle{\exp\left( -\frac\tau{2\eps^2}(\theta-\bar X_\tau)^2+U_\tau\right)}}{\displaystyle{\sum\limits_{k=1}^n \frac{e^{U_k} }{\sqrt k} }},
\end{equation}
where $U_k$ is defined in (\ref{Uk}). 
Chernoff and Zacks \cite{Chernoff&Zacks:1964} used Bayesian approach with normal priors on $\theta$ and uniform prior on $\tau$ to obtain an estimate of the mean of the observations after the change. In \cite{Lee&Heghiniam:1977} under the same assumptions the posterior distribution~(\ref{posterior}) in the change-point problem for normal observations was calculated.

Thus, we obtain the following Bayesian estimate of $L(\theta,\tau)$, 
\begin{equation}\label{Lb}
\widehat L_{\b}=\sum\limits_{\tau=1}^N \int\limits_\bR L(\theta,\tau)\; \pi(\theta,\tau|X)\,d\theta
= \frac1{\sqrt{2\pi\eps^2}} \sum_{\tau=1}^N p_\tau \sqrt\tau \int\limits_\bR L(\theta,\tau) \exp\left( -\frac\tau{2\eps^2}(\theta-\bar X_\tau)^2\right)\,d\theta,
\end{equation}
where
\begin{equation}\label{pk}
p_k=\frac{e^{U_k}}{\sqrt k}\left(\sum\limits_{i=1}^n \frac{e^{U_i} }{\sqrt i}\right)^{-1}.
\end{equation}
For example, if $L(\theta,\tau)=\theta\tau=\sum_{i=1}^\tau\theta_i$, then $\widehat L_\b$ is a weighted sum of $X_i$'s with weights $p_k$
$$
\widehat L_{\b}^\eps=\sum_{k=1}^n p_k \sum_{i=1}^k X_i\equiv \frac{\displaystyle{\sum\limits_{k=1}^n \frac{e^{U_k} }{\sqrt k} \sum\limits_{i=1}^ kX_i} }{\displaystyle{\sum\limits_{k=1}^n \frac{e^{U_k} }{\sqrt k}}}.
$$
Note that the Bayesian estimate of $\tau$ for quadratic loss function is given by
$$
\widehat\tau_{\mathrm {b}}^\eps=\sum_{k=1}^n k p_k\equiv\sum\limits_{k=1}^n \sqrt k e^{U_k}\left(\displaystyle{\sum\limits_{k=1}^n\frac{e^{U_k}}{\sqrt k} }\right)^{-1}
$$
and the corresponding estimate for $\theta$ is 
$$
\widehat\theta_\b^\eps=\sum_{k=1}^n p_k\bar X_k.
$$

\section{Simulation Study}

We studied the quadratic risks of Bayesian and maximum likelihood estimates of 
$$
L(\theta,\tau)=\sum_{i=1}^n\theta_i\equiv \theta\tau.
$$
$10^4$ simulations were made for $n=20$ observations in model~(\ref{model_chp}) with the values of $\theta\in\{0.5,1,1.5,2\}$, for the change-points $\tau=3,4,\dots,17,18$, and the noise level $\eps=1$.

\begin{figure}
\caption{Graphs of risk ratios $\kappa$ and $\tilde\kappa$ depending on $\tau\in\{3,4,\dots,18\}$  for $N=20$ observations, $\eps=1$, and different values of  $\theta=0.5,1,1.5,2$.}
\begin{center}
\subfigure[Graphs of risk ratio $\displaystyle{\kappa(\tau,\theta/\eps)=\frac{\E_\tau(\widehat\tau^\eps_{\mathrm{b}}-\tau)^2}{\E_\tau (\widehat\tau^\eps_{\mathrm{mle}}-\tau)^2}}$.] % caption for subfigure a
{
    \label{risk_graphs1:a}
    \includegraphics[width=12cm]{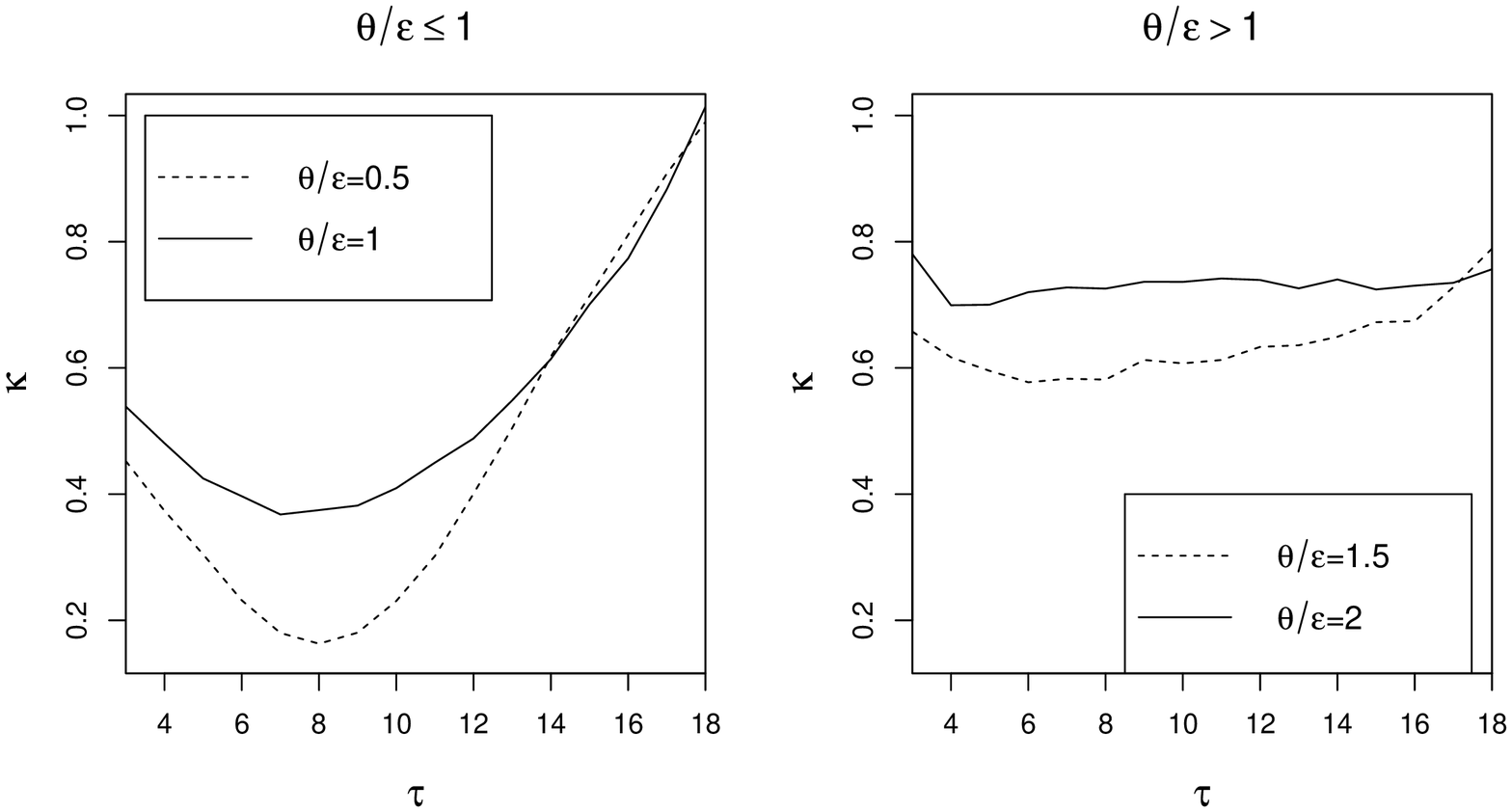}
}
\subfigure[Graphs of risk ratio $\displaystyle{\tilde\kappa(\tau,\theta/\eps)=\frac{\E_{\theta,\tau} (\widehat L^\eps_{\mathrm{b}}-L)^2}{\E_{\theta,\tau} (\widehat L^\eps_{\mathrm{mle}}-L)^2}}$.] % caption for subfigure b
{
    \label{risk_graphs1:b}
    \includegraphics[width=12cm]{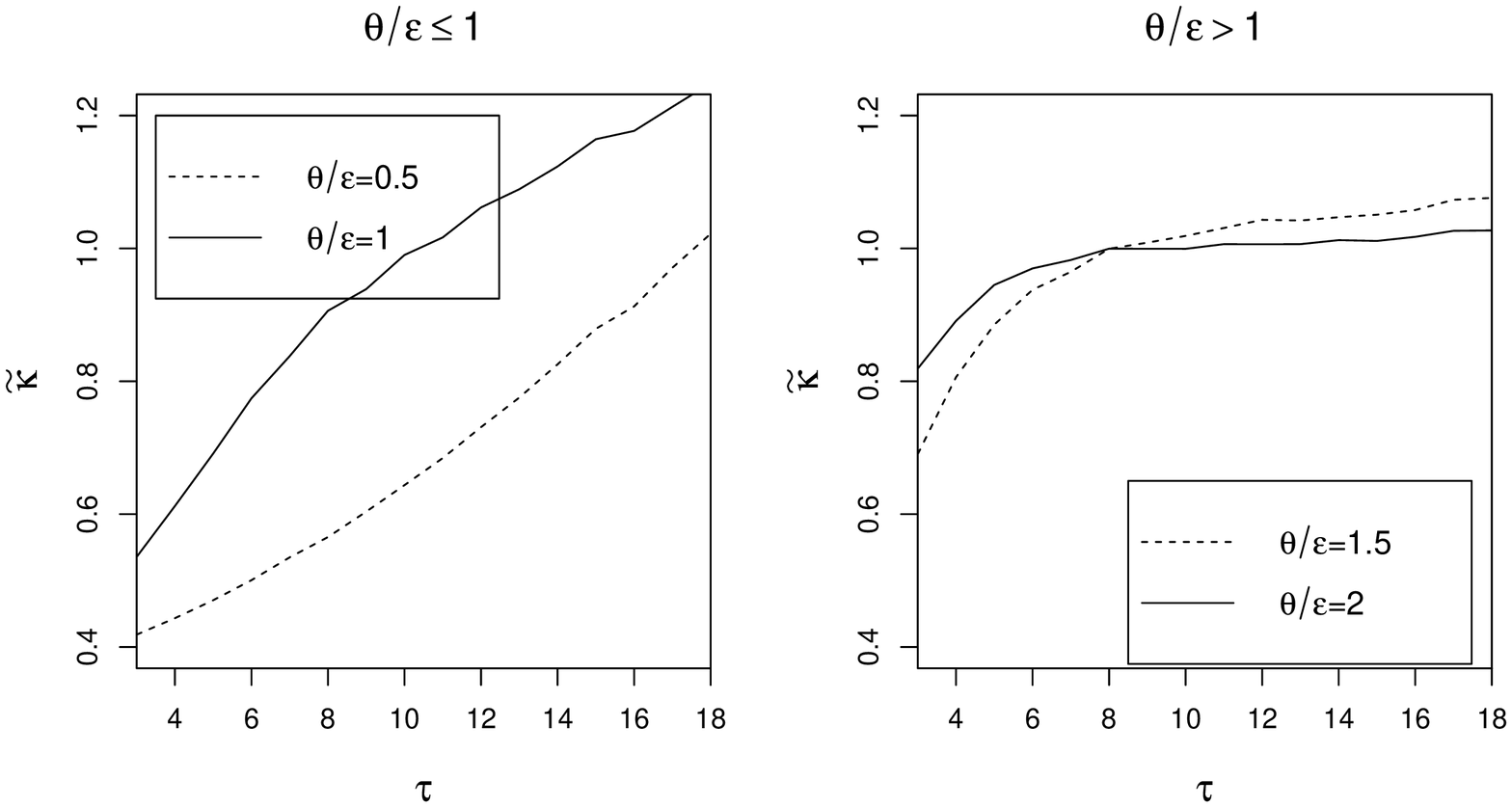}
}
\end{center}
\label{risk_graphs1} % caption for the whole figure
\end{figure}

First, introduce the following notation for risk ratios,
$$
\kappa(\tau,\theta/\eps)=\frac{\E_{\theta,\tau}  (\widehat \tau_{\b}^\eps-\tau)^2}{\E_{\theta,\tau}  (\widehat \tau_{\mle}^\eps-\tau)^2},\quad
\tilde \kappa(\tau,\theta/\eps)=\frac{\E_{\theta,\tau}  (\widehat L_{\b}^\eps-L)^2}{\E_{\theta,\tau}  (\widehat L_{\mle}^\eps-L)^2}.
$$
In Figure~\ref{risk_graphs1} the graphs of the empirical risk ratios $\kappa=\kappa(\tau,\theta/\eps)$ and $\tilde \kappa=\tilde\kappa(\tau,\theta/\eps)$  depending on $\tau=3,4,\dots,17,18$ are presented for different values of the signal-to-noise ratio (SNR), $\theta/\eps=0.5,1,1.5,2$. For simplicity we assume that the signal $\theta$ is positive.

Remind that asymptotically as $\eps\to 0$ the relative efficiency of the MLE of $\tau$ with respect to the Bayes estimate of $\tau$ is about 0.74,
$
\lim\limits_{\eps\to 0}\kappa(\tau,\theta/\eps)=\kappa_0\approx 0.7397
$.
It means that the MLE of $\tau$ is about $17\%$ less efficient than Bayesian estimate $\widehat\tau_{\mathrm{b}}$ if the SNR $\theta/\eps$ is large.

The examination of our numerical results leads to the following conclusions.
\begin{enumerate}
\item[(a)] {\bf Large SNR, $\theta/\eps>1$}.\\
It is clearly seen from Fig.~\ref{risk_graphs1}\subref{risk_graphs1:a}, that for large $\theta/\eps$ the ratio $\kappa(\tau,\theta/\eps)$ is close to its asymptotic theoretical value 0.7397. For $\theta/\eps=2$ the risk ratio $\kappa$ fluctuates between 0.72 and 0.77. For $\theta/\eps=1.5$ we have $0.57<\kappa<0.8$. However, this is not the case for the ratio $\tilde\kappa(\tau,\theta/\eps)$ of risks of estimating $L$ presented in Fig.~\ref{risk_graphs1}\subref{risk_graphs1:b}. In this case the behavior of the risk ratio depends both on $\theta/\eps$ and $\tau$.  For large SNR $\theta/\eps= 1.5,\ 2$ the relative efficiency is close to 1. It means that asymptotically both estimates of $L$ have very close risks. 
\item[(b)] {\bf Small SNR, $\theta/\eps\le 1$}.\\ 
For small SNR and moderate or small values of $\tau$, the Bayes estimate $\widehat L_{\mathrm{b}}$ has to be preferred to the MLE estimate $\widehat L_{\mathrm{mle}}$. For example, if $\theta/\eps=0.5$ and the change in the data takes place close to the beginning of the sequence  $\tau/N\le 0.4$, then 
 the Bayes estimate $\widehat L_{\mathrm{b}}$ of $L$ is almost twice more efficient than the MLE estimate $\widehat L_{\mathrm{mle}}$.   If $\tau$ is large, then the MLE of $L$ has to be chosen instead of the Bayes estimate. At the same time, the Bayes estimate of $\tau$ (Fig.~\ref{risk_graphs1}\subref{risk_graphs1:a}) is always more efficient than the MLE in the case of small SNR ($\theta/\eps=0.5,1$). Moreover, the smaller SNR is, the better is the behavior of Bayesian estimates comparing to the maximum likelihood estimates, both for estimating $\tau$ and $L$.
\item[(c)] {\bf Dependence on $\tau$}.\\ Fig.~\ref{risk_graphs1}\subref{risk_graphs1:b}  shows that the Bayesian estimate of $L$ is more efficient if the change takes place close to the beginning of the sequence. For example, for $\theta/\eps=1$ and $\tau/N<0.4$ the Bayes estimate is more efficient than MLE, and vice versa, the MLE of $L$ is more efficient for large values of $\tau$, $\tau/N>0.7$.  If the values of $\tau$ are moderate, $0.4\le \tau/N\le 0.7$, then depending on the SNR we should prefer MLE or the Bayes estimate of $L$ depending on the SNR.
\end{enumerate}
Our simulation results are very similar to the results of Sen and Srivastava~\cite{Sen&Srivastava:1975}. They made a comparative study of the Bayes and likelihood ratio tests for the problem of testing the hypothesis of "no change" in Gaussian data. It turned out that in the case of known mean $\theta$ in the data the Bayes test is superior for $\tau/N\le 0.4$, the LRT is superior for $\tau/N\ge 0.75$ and for $0.4<\tau/N<0.75$ the Bayes test dominates the LRT for small $\theta$ and vice versa.

The risks of MLE and Bayesian estimates of a smooth functional $\cL$ have the same first order asymptotic term as $\eps\to 0$. Thus, from the viewpoint of asymptotic behavior there is no difference what approach to choose for estimation. 

For small values of the signal-to-noise ratio Bayesian procedure has much better performance than the ML procedure for a large part of values of the change-point $\tau$ in case of quadratic losses.  We cannot explain this fact theoretically, since the behavior of the risk ratio is only known for large SNR as $\eps\to 0$.  

Simulation studies shows that Bayesian procedures work remarkably better than MLE procedures in the case of small signal-to-noise ratio. On the other hand, asymptotically, both procedures show the same performance. We think that due to this fact Bayesian estimates have to be used in non-asymptotic framework. Unfortunately, in non-asymptotic setting, their theoretical risk properties are very difficult to obtain.

{\bf Acknowledgements.}
The author is grateful to an anonymous referee for constructive comments and suggestions that helped to improve the paper.

\end{document}